\documentclass[12pt]{article}
\usepackage{amssymb}
\usepackage{amsfonts}
\usepackage{eurosym}
\usepackage{amsmath}

\newtheorem{theorem}{Theorem}

\newtheorem{corollary}[theorem]{Corollary}

\newtheorem{definition}[theorem]{Definition}

\newtheorem{lemma}[theorem]{Lemma}

\newtheorem{proposition}[theorem]{Proposition}

\newenvironment{proof}[1][Proof]{\noindent\textbf{#1.} }{\ \rule{0.5em}{0.5em}}
\input{tcilatex}
\begin{document}

\setlength{\baselineskip}{1.5\baselineskip}\newpage

\begin{center}
{\large On Fixed Points of Locally and Pointwise Contractive Set-Valued Maps
with an Application to the Existence of Nash Equilibrium in Games\footnote{%
I thank the Editor, participants of Rice University's Computational
Mathematics and Operations Research Colloquium Series, as well as Dick
McGehee, Vladimir Sverak, and especially Peter Olver for comments and
suggestions on earlier drafts. Any remaining errors are mine.}}

\bigskip \bigskip \bigskip \bigskip

Ted Loch-Temzelides

{\small Rice University}

{\small 6100 Main Street, Houston, TX 77005, USA}

{\small tedt@rice.edu}

\bigskip \bigskip \bigskip \bigskip \bigskip \bigskip \bigskip

\textbf{Abstract}
\end{center}

\noindent We establish the existence of fixed points for set-valued maps
defined on metric spaces and satisfying a pointwise or a local version of
Banach's contraction property. As an application, we demonstrate the
existence of Nash equilibrium in a general class of strategic games played
on metric spaces of strategies.

\begin{center}
\bigskip \bigskip \bigskip
\end{center}

\newpage

\section{Introduction}

The study of fixed points in metric spaces constitutes an important part of
analysis and has many uses within and outside of mathematics. The modern
theory of games originates in von Neumann and Morgenstern (1947). Strategic
games have several applications of game theory to biology, control theory,
engineering, economics, and related fields. The existence of Nash
equilibrium for finite games, viewed as fixed points of certain set-valued
maps, was first demonstrated by Nash (1950). Related proofs often use some
version of the fixed-point theorem in Kakutani (1941). For games with
infinite-dimensional strategy spaces, existence results often employ a
related generalization by Glicksberg (1952). In addition to requiring that
the set-valued map is convex-valued and has a closed graph, these results
rely on the assumption that the domain of the set-valued map is convex.

We extend recent results on the existence of fixed points for maps
satisfying pointwise and local contraction properties to set-valued maps
defined on metric spaces. As an illustration, we use these results to
demonstrate the existence of pure-strategy Nash equilibrium for a certain
class of strategic games, which we term \textit{contractive games}. Our
treatment allows for the domain of the set-valued map to be rather general
(possibly infinitely dimensional) metric spaces. Our most general results
require that the set-valued map is compact-valued, that it satisfies a
pointwise version of Banach's contraction mapping property, and that its
domain is a compact, rectifiable path connected metric space.

Local and pointwise contraction properties are, of course, weaker than
assuming that the set-valued map is a contraction globally.\ As a result,
completeness of the underlying metric space will not be sufficient for the
existence of a fixed point, thus, that of a Nash equilibrium. While we use
existence of Nash equilibrium in strategic games as an illustration of our
main results, these techniques can be used in other applications, including
the existence of value functions in set-valued dynamic programming.

\section{Review of Main Concepts}

A classic result by Banach (1922) established the existence of a unique
fixed point in complete metric spaces for maps satisfying a contraction
property. Since the beginning of the twentieth century, there has been a
well-developed theory of set-valued maps, that is maps $F:X\rightarrow 2^{X}$%
, see, for example, Kuratowski (1958). Selected recent references, to name a
few, include Gorniewicz (2006), Aubin and Frankowska (2009), Arutyunov and
Obukhovskii (2017), and Malkoun and Olver (2021). Fixed points for
set-valued maps are defined as points $x^{\ast }$ such that $x^{\ast }\in
F(x^{\ast })$. A celebrated result due to Kakutani (1941) established the
existence of fixed points for set-valued maps satisfying certain continuity
and convexity properties. Kakutani's result relies on the domain being
convex (a finite-dimensional simplex). The notion of a contraction can be
applied to set-valued maps, for example by using a Hausdorff distance. Using
this notion, generalizations of Banach's existence result for set-valued
maps can be obtained as in Nadler (1969). Unlike in the single-valued case,
uniqueness is not guaranteed.

There have been attempts to generalize Banach's result for single-valued
maps satisfying weaker, pointwise or local notions of the standard
contraction property; see Ciesielski and Jasinski (2018) for a recent
comprehensive review. These results typically require additional conditions
to the completeness of the underlying metric space, such as compactness and
some notion of connectedness. These restrictions are topological in nature
and do not require convexity. Here we investigate under what conditions
pointwise or local contraction properties are sufficient for the existence
of fixed points for set-valued maps. The common thread in the argument used
involves two steps. First, using a Hausdorff distance provides us with a
complete metric space. We can then define set-valued maps and their
pointwise or local contraction properties in that space. The second step
involves imposing additional topological properties on the space and
defining a new metric under which the underlying metric space remains
complete, and the local or pointwise set-valued contraction becomes a global
set-valued contraction. Related properties have been previously established
for maps, and this paper demonstrates that, subject to appropriate
modifications, they can be extended to set-valued maps.

Throughout the paper, we let $(X,d)$, $X\neq \varnothing $, be a metric
space. Additional restrictions on $X$ will be imposed as needed in what
follows. Banach (1922) demonstrated that if $(X,d)$ is complete, and $%
f:X\rightarrow X$ satisfies the \textit{contraction property} that there is $%
\beta \in \lbrack 0,1)$ such that for all $x,y\in X$, $d(f(x),f(y))\leq
\beta d(x,y)$, then there exists a unique $x^{\ast }\in X$ such that $%
x^{\ast }=f(x^{\ast })$. The map $f$ is sometimes referred to as a \textit{%
contraction of modulus} $\beta $. A number of generalizations of Banach's
contraction property have been investigated. In what follows, we will
concern ourselves with the implications of extending the local and pointwise
versions studied in Ciesielski and Jasinski (2018) to the study of
set-valued maps. We will make use of the following notions.

\begin{definition}
$X$ is \textbf{convex} if for all $x,y\in X$, there exists $z\in X$ such
that $d(x,y)=d(x,z)+d(z,y)$.
\end{definition}

In the case where $X$ is complete and convex, whenever $x,y\in X$, then $X$
also contains a subset whose boundary points are $x$ and $y$, and which is
isometric to an interval of length $d(x,y)$. The following notion of
convexity is commonly employed in metric spaces.

\begin{definition}
$X$ is $d$\textbf{-convex} if for all $x,y\in X$, there exists a
homeomorphism (path) $p:[0,1]\rightarrow X$ with $p(0)=x$, $p(1)=y$ and such
that for any $0\leq t_{1}\leq t_{2}\leq t_{3}\leq 1$ we have $%
d(p(t_{1}),p(t_{3}))=d(p(t_{1}),p(t_{2}))+d(p(t_{2}),p(t_{3}))$.
\end{definition}

Thus, if $X$ is $d$-convex then for any $x,y\in X$, there is an interval $%
[x,y]\subseteq X$ that is the image of a homeomorphism $p:[0,1]\rightarrow X$%
.

\begin{definition}
Fix $r>0$. $X$ is $r$\textbf{-chainable} if given any $x,y\in X$, there are
finitely many points $z_{0},\ldots z_{n}\in X$ such that $z_{0}=x$, $z_{n}=y$%
, and $d(z_{i-1},z_{i})<r$, for all $i=1,\ldots ,n$.
\end{definition}

Finally, we will employ the following.

\begin{definition}
Let $\Pi $ be the set of all possible partitions $%
0=t_{0}<t_{1}<...<t_{n-1}<t_{n}=1$ of $[0,1]$, $n<\infty $. $X$ is \textbf{%
(rectifiably) path-connected} if for all $x,y\in X$, there exists a
continuous function (path) $p:[0,1]\rightarrow X$ such that $p(0)=x$, $p(1)=y
$, and whose length 
\begin{equation*}
l(p):=\sup_{\pi \in \Pi }\dsum\limits_{i=1}^{n}d(p(t_{i-1}),p(t_{i}))<\infty 
\end{equation*}
\end{definition}

In other words, $X$ is rectifiably path-connected if every two points of $X$
can be connected by a path $p$ in $X$ of finite length, $l(p)$. Later we
will use $l(p^{[a,b]})$ to indicate the restriction of a path $%
p:[0,1]\rightarrow X$ on $[a,b]$, where $0\leq a<b\leq 1$.

While Banach's Theorem makes no use of connectedness properties, the
situation is more delicate once only local or pointwise contraction
properties are assumed. Indeed, absent the rectifiable path connectedness
assumption, Ciesielski and Jasinskic (2016a) used a Cantor set-like
structure to construct a compact $X\subset 
%TCIMACRO{\U{211d} }%
%BeginExpansion
\mathbb{R}
%EndExpansion
$ containing no isolated points, and a smooth bijection $f:X\rightarrow X$
such that $f$ has zero slope everywhere, but it admits no fixed (or
periodic) points.

In what follows, we will restrict attention to set-valued maps. Let $(X,d)$
be a metric space with $X\neq \varnothing $. A \textit{set-valued map} is a
map $F:X\rightarrow 2^{X}\backslash \varnothing $. A \textit{fixed point of }%
$F$ is an $x^{\ast }\in X$ such that $x^{\ast }\in F(x^{\ast })$. We will
let \emph{Fix}$(F)$ stand for the set of fixed points of $F$ in $X$. We
begin with some basic concepts related to set-valued maps.

Let $C(X):=\{Y\subseteq X:Y\neq \varnothing $, $Y$ closed and bounded$\}$,
and let $K(X):=\{Y\subseteq X:Y\neq \varnothing $, $Y$ compact$\}$. Let $d$
be a metric on $X$. For any $X_{1},X_{2}\in $ $C(X)$ the \textit{Hausdorff
distance} is defined as follows. Let

\begin{equation}
d(y,X_{1}):=\inf_{x\in X_{1}}d(x,y)
\end{equation}%
Then:

\begin{equation}
H(X_{1},X_{2}):=\max \left\{ \sup_{x\in X_{1}}d(x,X_{2}),\sup_{y\in
X_{2}}d(y,X_{1})\right\}
\end{equation}%
Equivalently, letting

\begin{equation}
B_{r}(X_{1}):=\dbigcup\limits_{x\in X_{1}}\left\{ z\in X:d(x,z)\leq r\right\}
\end{equation}%
we have:

\begin{equation}
H(X_{1},X_{2}):=\inf \left\{ r:X_{1}\subset B_{r}(X_{2})\text{ and }%
X_{2}\subset B_{r}(X_{1})\right\} 
\end{equation}%
Using the Hausdorff distance, Nadler (1969) generalized the notions of
Lipschitz and contraction maps to set-valued maps.

\begin{definition}
Let $(X,d_{X})$ and $(Y,d_{Y})$ be metric spaces. A set-valued map $%
F:X\rightarrow Y$ is a \textbf{set-valued Lipschitz map}\textit{\ }if for
all $x,y\in X$, $H(F(x),F(y))\leq \alpha d(x,y)$, where $\alpha \in \lbrack
0,\infty )$ is the \textbf{Lipschitz constant}. $F$ is a \textbf{set-valued
contraction of modulus }$\beta $\textit{\ }if it is a set-valued Lipschitz
map with Lipschitz constant $\beta \in \lbrack 0,1)$.
\end{definition}

Nadler (1969) demonstrated the existence of fixed points for set-valued maps
defined on complete metric spaces and satisfying the above contraction
property. Several generalizations have been studied since; see, for example,
Reich (1972a,b). We now turn to the main focus of this paper, which concerns
extending results on the existence of fixed points of pointwise and local
contractions to set-valued maps. The following defines iterates of
set-valued maps. Let $X$, $Y$, $Z$ be metric spaces.

\begin{definition}
Let $F_{1}:X\rightarrow Y$ and $F_{2}:Y\rightarrow Z$. The \textbf{%
composition} of $F_{1}$ and $F_{2}$ is defined by: $(F_{2}\circ
F_{1})(x):=\cup \{F_{2}(y):y\in F_{1}(x)\}$, for all $x\in X$. Proceeding
recursively, we obtain: $F^{(n)}(x):=$ $\underset{n\text{-times}}{%
\underbrace{F\circ F\ldots \circ F}}$.
\end{definition}

We remark that if $F_{1}$ is a contraction of modulus $\beta _{1}$ and $%
F_{2} $ is a contraction of modulus $\beta _{2}$, then $(F_{2}\circ F_{1})$
is a contraction map of modulus $\beta _{1}\cdot \beta _{2}$.

The following definitions extend the related notions in Ciesielski and
Jasinski (2018) to set-valued maps.

\begin{definition}
A set-valued map $F:X\rightarrow X$ is a \textbf{set-valued pointwise
contraction}\textit{\ }if for all $x\in X$ there exists $\beta _{x}\in
\lbrack 0,1)$ and an open neighborhood $N(x)\subseteq X$ such that for all $%
y\in N(x)$, $H(F(x),F(y))\leq \beta _{x}d(x,y)$. The set-valued map $F$ is a 
\textbf{set-valued} \textbf{uniform pointwise contraction }if it is a
set-valued pointwise contraction and $\beta $ does not depend on $x$.
\end{definition}

A stronger property is defined in the following.

\begin{definition}
A set-valued map $F:X\rightarrow X$ is a $(\beta _{x},r_{x})$-\textbf{local
contraction}\textit{\ }if for all $x\in X$ there exists $r_{x}>0$ and $\beta
_{x}\in \lbrack 0,1)$ such that $F$ is a set-valued contraction with modulus 
$\beta _{x}$ in $B_{r_{x}}(x)$. The set-valued map $F$ is a $(\beta ,r)$-%
\textbf{uniform local contraction} if $\beta $ and $r$ do not depend on $x$.
\end{definition}

A pointwise/local contraction is also a uniform pointwise/uniform local
contraction if the underlying space is compact. A (uniform) local
contraction is a (uniform) pointwise contraction, but the reverse
implication does not necessarily hold.

\section{Fixed Points of Locally Contractive Set-Valued Maps}

As before, we let $(X,d)$ be a metric space with $X\neq \varnothing $. We
have the following.

\begin{proposition}
Let $X$ be convex and complete and let $F:X\rightarrow K(X)$ be a set-valued 
$(\beta ,r)$-uniform local contraction on $X$. Then \emph{Fix}$(F)\neq
\varnothing $.
\end{proposition}

\begin{proof}
It is sufficient to show that $F$ is a set-valued contraction of modulus $%
\beta $ on $X$. The conclusion then follows from Nadler (1969). To this end,
let $x,y\in X$. By convexity of $X$, given any $r>0$, there are points $%
z_{0}=x,z_{1},\ldots ,z_{n-1},z_{n}=y$ such that $d(x,y)=%
\sum_{i}d(z_{i-1},z_{i})$ and $d(z_{i-1},z_{i})<r$. Then%
\begin{eqnarray}
H\left( F(x),F(y)\right) &\leq &\sum_{i=1}^{n}H\left(
F(z_{i-1}),F(z_{i})\right)  \notag \\
&<&\beta \sum_{i=1}^{n}d(z_{i-1},z_{i})  \notag \\
&=&\beta d(x,y)
\end{eqnarray}%
Thus, $F$ is a set-valued contraction.
\end{proof}

Next, we assume that $X$ is compact and $r$-chainable and present an
existence proof for the case where $F$ is a (not necessarily uniform)
locally contractive set-valued map. We will use the following.

\begin{lemma}
Let $X$ be $r$-chainable. For each $x$, $y\in X$, define the function $%
d_{c}:X\times X\rightarrow 
%TCIMACRO{\U{211d} }%
%BeginExpansion
\mathbb{R}
%EndExpansion
_{+}$ by%
\begin{equation}
d_{c}(x,y):=\inf \left\{ \sum_{i=1}^{n}d(z_{i-1},z_{i}):z_{0}=x,z_{1},\ldots
,z_{n}=y\text{ is an }r\text{-chain from }x\text{ to }y\right\}
\end{equation}%
Then: (a) $d_{c}:X\times X\rightarrow 
%TCIMACRO{\U{211d} }%
%BeginExpansion
\mathbb{R}
%EndExpansion
_{+}$ is a metric on $X$ and (b) $d_{c}$ is topologically equivalent to $d$.
Thus, the space $(X,d_{c})$ is complete whenever $(X,d)$ is complete.
\end{lemma}

\begin{proof}
\textit{(a)} Clearly, $d_{c}(x,y)\geq 0$, $d_{c}(x,y)=0$ if and only if $x=y$%
, and $d_{c}(x,y)=d_{c}(y,x)$ by construction. For the triangle inequality,
fix $x,y,z\in X$ and $\epsilon >0$. Let $w_{0}=x,w_{1},\ldots
,w_{n-1},w_{n}=y$ be an $r$-chain from $x$ to $y$ and $v_{0}=y,v_{1},\ldots
,v_{n-1},v_{n}=z$ be an $r$-chain from $y$ to $z$ such that 
\begin{eqnarray}
d_{c}(x,y) &\geq &\sum_{i=1}^{n}d(w_{i-1},w_{i})-\epsilon \\
d_{c}(y,z) &\geq &\sum_{i=1}^{n}d(v_{i-1},v_{i})-\epsilon
\end{eqnarray}%
Then $w_{0},w_{1},\ldots ,w_{n-1},w_{n},v_{0},v_{1},\ldots ,v_{n-1},v_{n}$
is an $r$-chain from $x$ to $z$ and 
\begin{eqnarray}
&&d_{c}(x,y)+d_{c}(y,z)  \notag \\
&\geq &\sum_{i=1}^{n}d(w_{i-1},w_{i})-\epsilon
+\sum_{i=1}^{n}d(v_{i-1},v_{i})-\epsilon  \notag \\
&\geq &d_{c}(x,z)-2\epsilon
\end{eqnarray}%
Since $\epsilon >0$ was arbitrary, the triangle inequality follows.

\textit{(b)} Clearly, $d(x,y)\leq d_{c}(x,y)$, for all $x$, $y\in X$. In
addition, when $d(x,y)<r$, we have 
\begin{gather}
d(x,y)\leq d_{c}(x,y)  \notag \\
=\inf \left\{ \sum_{i=1}^{n}d(z_{i-1},z_{i}):z_{0}=x,\ldots ,z_{n}=y\text{
is an }r\text{-chain from }x\text{ to }y\right\}  \notag \\
=\inf \left\{ d(z_{0},z_{1}):z_{0}=x,z_{1}=y\right\}  \notag \\
=d(x,y)
\end{gather}%
Hence, the two metrics are topologically equivalent and $(X,d_{c})$ is a
complete metric space whenever $(X,d)$ is complete.
\end{proof}

\bigskip

Let $X$ be compact and $F:X\rightarrow K(X)$ be a set-valued $(\beta
_{x},r_{x})$-local contraction on $X$. The set $\mathcal{\{}%
B_{r_{x}}(x):x\in X\}$ forms an open covering of $X$. Since $X$ is compact,
it has a finite subcover: $\mathcal{B}=\{B_{r_{x}}^{\prime }(x):x\in X\}$.
Let $\beta =\max \{\beta _{x}:B_{r_{x}}^{\prime }(x)\in \mathcal{B}\}<1$ and 
$0<r=\min \left\{ r_{x}:x\in X\right\} $. Then $F$ is a $(\beta ,r)$-uniform
local contraction. Given such an $F$, we have the following.

\begin{proposition}
Let $X$ be compact and let $F:X\rightarrow K(X)$ be a set-valued $(\beta
_{x},r_{x})$-local contraction on $X$. Then: (a) if $X$ is $r$-chainable,
then \emph{Fix}$(F)\neq \varnothing $; (b) if $X$ is the finite union of
connected components that are pairwise disjoint and $F(x)$ is connected, for
all $x$, then \emph{Fix}$(F^{(n)})\neq \varnothing $.
\end{proposition}

\begin{proof}
\textit{(a)} By Lemma 10, it suffices to show that the local contraction
property implies that $F:X\rightarrow K(X)$ is a contraction with respect to
metric $d_{c}$ and the corresponding Hausdorff distance, $H_{c}$. Let $x$, $%
y\in X$ and $z_{0}=x,z_{1},\ldots ,z_{n-1},z_{n}=y$ be an $r$-chain from $x$
to $y$. Since $d(z_{i-1},z_{i})<r$, for all $i$, we have that $%
H(F(z_{i-1}),F(z_{i}))\leq \beta d(z_{i-1},z_{i})<r$, where $r$ and $\beta $
are as discussed above. Thus,%
\begin{eqnarray}
H_{c}(F(x),F(y)) &\leq &\sum_{i=1}^{n}H_{c}(F(z_{i-1}),F(z_{i}))  \notag \\
&=&\sum_{i=1}^{n}H(F(z_{i-1}),F(z_{i}))  \notag \\
&\leq &\beta \sum_{i=1}^{n}d(z_{i-1},z_{i})  \notag \\
&=&\beta d_{c}(x,y)
\end{eqnarray}%
Since $z_{0}=x,z_{1},\ldots ,z_{n-1},z_{n}=y$ was an arbitrary $r$-chain, we
have that $H_{c}(F(x),F(y))\leq \beta d_{c}(x,y)$; i.e., $F$ is a set-valued
contraction with respect to $d_{c}$ and $H_{c}$.

\textit{(b) }Now suppose that $X=X_{1}\cup ...\cup X_{n}$, where $X_{i}$ is
connected, for all $i$, and $X_{i}\cap X_{j}=\varnothing $, $i\neq j$. Since 
$F$ is a contraction and $F(x)$ is connected, we have that $F^{(n)}(X_{j})$
is connected, for any $j$. Since $X$ consists of finitely many disconnected
components, for all $j$, there exists $k$, $l$, and $m$, such that $%
F^{(k)}(X_{j})\cap X_{m}\neq \varnothing $ and $F^{(k+l)}(X_{j})\cap
X_{m}\neq \varnothing $, where $X_{m}$ is the same connected component of $X$%
. Thus, $F^{(l)}(X_{m})\subset X_{m}$. Since $X_{m}$ is connected, it is $r$%
-chainable. Applying part \textit{(a)} of the Proposition to $F^{(l)}$ on $%
X_{m}$, we conclude that there exists an $x^{\ast }\in X_{m}\subset X$ such
that $x^{\ast }\in F^{(l)}(x^{\ast })$.
\end{proof}

\bigskip

\section{Fixed Points of Pointwise Contractive Set-Valued Maps}

In what follows we consider two cases. First, we will assume that $(X,d)$ is 
$d$-convex. Then we will impose the condition that $X$ is rectifiably
path-connected. For the first case, we have the following.

\begin{proposition}
Let $X$ be complete and $d$-convex and let $F:X\rightarrow K(X)$ be a
uniform pointwise set-valued contraction on $X$. Then \emph{Fix}$(F)\neq
\varnothing $.
\end{proposition}

\begin{proof}
Since $F$ is a uniform pointwise contraction, for all $x\in X$ there exists $%
\beta \in \lbrack 0,1)$ and an open neighborhood $N_{\epsilon }(x)\subseteq
X $ such that if $y\in N_{\epsilon }(x)$, then $H(F(x),F(y))\leq \beta
d(x,y) $. It is sufficient to show that for all $x\in X$ there exists $\beta
\in \lbrack 0,1)$ such that $F$ is a set-valued contraction with modulus $%
\beta $ on $X$.

Let $x,y\in X$. Since $X$ is $d$-convex, there exists a continuous $%
p:[0,1]\rightarrow X$ with $p(0)=x$ and $p(1)=y$. Since $F$ is a uniform
pointwise contraction on $X$, for each $t\in \lbrack 0,1]$ there is $%
\epsilon _{t}>0$ such that if $d(p(t),x)<\epsilon _{t}$, then $H(F\circ
p(t),F(x))\leq \beta d(p(t),x)$. Since $p$ is continuous, for all $t\in
(0,1) $, there exists $\delta _{t}>0$ with $B_{\delta _{t}}(t)\subseteq
(0,1) $ and such that $t^{\prime }\in B_{\delta _{t}}(t)$ implies $%
d(p(t),p(t^{\prime }))<\epsilon _{t}$. Similarly, for the endpoints $t=0,1$
we can choose $\delta _{0},\delta _{1}>0$ such that $B_{\delta
_{0}}(0)\subseteq \lbrack 0,1)$, $B_{\delta _{1}}(1)\subseteq (0,1]$, and $%
t^{\prime }\in B_{\delta _{i}}(i)$ implies $d(p(i),p(t^{\prime }))<\epsilon
_{i}$, $i=0,1$. The collection $\{B_{\delta _{t}}(t)\subset (0,1):t\in
\lbrack 0,1]\}$ forms an open cover of $[0,1]$. Since $[0,1]$ is connected,
there is a sequence of sets $\{B_{\delta _{t_{i}}}(t_{i})\}$ such that $%
[0,1]\subset \cup _{i=1}^{n}B_{\delta _{t_{i}}}(t_{i})$, $0\in B_{\delta
_{0}}(t_{0})$, $1\in B_{\delta _{1}}(t_{1})$, and $B_{\delta
_{t_{i}}}(t_{i})\cap B_{\delta _{t_{j}}}(t_{j})\neq \varnothing $, if $%
|i-j|\leq 1$.

Since $B_{\delta _{t_{i}}}(t_{i})$ is a symmetric interval centered at $t_{i}
$, we can choose $c_{i}\in B_{\delta _{t_{i-1}}}(t_{i-1})\cap B_{\delta
_{t_{i}}}(t_{i})$ so that $t_{i-1}<c_{i}<t_{i}$, $i=1,...,n$. Choosing $%
\delta _{t_{i}}>0$ such that $d(p(t_{i}),p(t_{i}^{\prime }))<\epsilon
_{t_{i}}$ if $t_{i}^{\prime }\in B_{\delta _{t_{i}}}(t_{i})$, we have that $%
d(p(t_{i-1}),p(c_{i}))<\epsilon _{t_{i-1}}$ and $d(p(c_{i}),p(t_{i}))<%
\epsilon _{t_{i}}$. Since $F\ $is a uniform pointwise contraction, we have
that, for all $i$,%
\begin{eqnarray}
H(F\circ p(t_{i-1}),F\circ p(t_{i})) &\leq &H(F\circ p(t_{i-1}),F\circ
p(c_{i}))+H(F\circ p(c_{i}),F\circ p(t_{i}))  \notag \\
&\leq &\beta d(p(t_{i-1}),p(c_{i}))+\beta d(p(c_{i}),p(t_{i}))  \notag \\
&=&\beta d(p(t_{i-1}),p(t_{i}))
\end{eqnarray}%
Thus,%
\begin{eqnarray}
H(F(x),F(y)) &=&H(F\circ p(t_{0}),F\circ p(t_{n}))  \notag \\
&\leq &\sum_{i=1}^{n}H(F\circ p(t_{i-1}),F\circ p(t_{i}))  \notag \\
&\leq &\beta \sum_{i=1}^{n}d(p(t_{i-1}),p(t_{i}))  \notag \\
&\leq &\beta d(p(0),p(1))=\beta d(x,y)
\end{eqnarray}%
Hence, $F$ is a set-valued contraction on $X$.
\end{proof}

\bigskip

Next, we assume that $X$ is rectifiably path-connected. Recall that, given $%
0<a<b<1$, a rectifiable path is a continuous map from a closed real interval 
$[a,b]$ into $X$, whose length is given by 
\begin{equation}
l(p^{[a,b]})=\sup_{\pi \in \Pi }\left\{
\dsum\limits_{i=1}^{n}d(p(t_{i-1}),p(t_{i})):\pi
=\{t_{0}=a,...,t_{n}=b\}\right\} <\infty
\end{equation}%
If $p$ is a rectifiable path, then for all $c\in \lbrack a,b]$ we have that $%
l(p^{[a,b]})=l(p^{[a,c]})+l(p^{[c,b]})$. First, using the Hausdorff
distance, we extend the notion of length to set-valued paths in order to
demonstrate the following.

\begin{lemma}
Let $F:X\rightarrow K(X)$ be a set-valued uniform pointwise contraction with
modulus $\beta \in \lbrack 0,1)$ and let $p:[0,1]\rightarrow X$ be a path
with $l(p)<\infty $. Then $l(F\circ p)\leq \beta l(p)<\infty $, where 
\begin{equation*}
l(F\circ p):=\sup_{\pi \in \Pi }\dsum\limits_{i=1}^{n}H(F\circ
p(t_{i-1}),F\circ p(t_{i}))
\end{equation*}
\end{lemma}

\begin{proof}
Let $F:X\rightarrow K(X)$ be a set-valued uniform pointwise contraction and
consider a rectifiable path $p:[0,1]\rightarrow X$. We first let $0<a<b<1$
and show that $H(F\circ p(a),F\circ p(b))\leq \beta l(p^{[a,b]})$.
Rectifiability and the contraction property imply a partition $\pi
=\{t_{i}:t_{0}=a,...t_{i-1},t_{i},...t_{n}=b\}$ and $s_{i}\in (t_{i-1},t_{i})
$, $i=1,...,n$, such that%
\begin{eqnarray}
H(F\circ p(t_{i-1}),F\circ p(t_{i})) &\leq &H(F\circ p(t_{i-1}),F\circ
p(s_{i}))+H(F\circ p(s_{i}),F\circ p(t_{i}))  \notag \\
&\leq &\beta \lbrack d(p(t_{i-1}),p(s_{i}))+d(p(s_{i}),p(t_{i}))]
\end{eqnarray}%
For all $0<s_{1}<s_{2}<1$, $d(s_{1},s_{2})\leq l(p^{[s_{1},s_{2}]})$. Thus,
(15) implies that, for all $i=1,...,n$, 
\begin{equation}
H(F\circ p(t_{i-1}),F\circ p(t_{i}))\leq \beta \lbrack
l(p^{[t_{i-1},s_{i}]})+l(p^{[s_{i},t_{i}]})]=\beta l(p^{[t_{i-1},t_{i}]})
\end{equation}%
Hence,%
\begin{gather}
H(F\circ p(a),F\circ p(b))=H(F\circ p(t_{0}),F\circ p(t_{n}))  \notag \\
\leq \dsum\limits_{i=1}^{n}H(F\circ p(t_{i-1}),F\circ p(t_{i}))  \notag \\
\leq \beta \dsum\limits_{i=1}^{n}l(p^{[t_{i-1},t_{i}]})=\beta l(p^{[a,b]})
\end{gather}%
Thus, $H(F\circ p(a),F\circ p(b))\leq \beta l(p^{[a,b]})$, for all $0<a<b<1$%
. Next, consider $\pi =\{t_{i}:t_{0}=0,...t_{i-1},t_{i},...t_{n}=1\}$.
Expression (16) gives%
\begin{equation}
\dsum\limits_{i=1}^{n}H(F\circ p(t_{i-1}),F\circ p(t_{i}))\leq \beta
\dsum\limits_{i=1}^{n}l(p^{[t_{i-1},t_{i}]})=\beta l(p)
\end{equation}%
Hence, $\beta l(p)$ is an upper bound of the set 
\begin{equation}
\left\{ \dsum\limits_{i=1}^{n}H(F\circ p(t_{i-1}),F\circ p(t_{i})):\pi
=\{t_{0}=0,...,t_{n}=1\}\right\} 
\end{equation}%
By the least upper bound property of $%
%TCIMACRO{\U{211d} }%
%BeginExpansion
\mathbb{R}
%EndExpansion
$, the $\sup $ of this set, $l(F\circ p)$, exists and satisfies $l(F\circ
p)\leq \beta l(p)$.
\end{proof}

\bigskip

We now have the following.

\begin{proposition}
Let $X$ be complete and rectifiably path-connected and let $F:X\rightarrow
K(X)$ be a set-valued uniform pointwise contraction with modulus $\beta \in
\lbrack 0,1)$. Then \emph{Fix}$(F)\neq \varnothing $.
\end{proposition}

\begin{proof}
Let $x$, $y\in X$ and, as before and let $p:[0,1]\rightarrow X$ be such that 
$p(0)=x$ and $p(1)=y$. Define:%
\begin{gather}
d_{r}(x,y):=\inf_{p}\left\{ l(p):p\text{ is a rectifiable path from }x\text{
to }y\right\}  \notag \\
=\inf_{p}\sup_{\pi \in \Pi
}\dsum\limits_{i=1}^{n}d(p(t_{i-1}),p(t_{i}))<\infty
\end{gather}%
The function $d_{r}:X\times X\rightarrow 
%TCIMACRO{\U{211d} }%
%BeginExpansion
\mathbb{R}
%EndExpansion
_{+}$ is a metric on $X$. Let $H_{r}$ be the Hausdorff distance associated
with $d_{r}$. It remains to show that: \textit{(a)} $F:X\rightarrow K(X)$ is
a contraction with respect to the metrics $d_{r}$ and $H_{r}$, and \textit{%
(b)} $(X,d_{r})$ is complete, if $(X,d)$ is complete. As before, the
existence of a fixed point of $F$ then follows from Nadler (1969).

For \textit{(a)}, let $x,y\in X$, and $\epsilon >0$ be given. Since $X$ is
rectifiably path-connected, there exists a rectifiable path $p$ from $x$ to $%
y$ such that $l(p)\leq d_{r}(x,y)+\epsilon $. Then, (20) and Lemma 13 imply%
\begin{eqnarray}
H_{r}(F(x),F(y)) &\leq &\sup_{\pi \in \Pi }\dsum\limits_{i=1}^{n}H(F\circ
p(t_{i-1}),F\circ p(t_{i}))  \notag \\
&=&l(F\circ p)\leq \beta l(p)\leq \beta (d_{r}(x,y)+\epsilon )
\end{eqnarray}%
Since $\epsilon $ was arbitrary, $H_{r}(F(x),F(y))\leq \beta d_{r}(x,y)$,
for all $x,y\in X$. Thus, $F$ is a contraction on $X$.

For \textit{(b)}, let $\{x_{n}\}$ be a Cauchy sequence in $(X,d_{r})$. Since 
$d(x,y)\leq d_{r}(x,y)$, for all $x,y\in X$, $\{x_{n}\}$ is also a Cauchy
sequence in $(X,d)$. Since $(X,d)$ is complete, there exists $\overline{x}%
\in X$ such that $\lim_{n\rightarrow \infty }d(x_{n},\overline{x})=0$. It is
sufficient to show that $\lim_{n\rightarrow \infty }d_{r}(x_{n},\overline{x}%
)=0$. Since $\{x_{n}\}$ is Cauchy in $(X,d_{r})$, there exists a sequence $%
\{\epsilon _{i}\}$, where $\epsilon _{i}>0$, for all $i$, and $\sum
{}_{i=1}^{\infty }\epsilon _{i}<\infty $, and an increasing sequence $%
\{N_{i}\}$ such that if $m,n\geq N_{i}$ then $d_{r}(x_{n},x_{m})<\epsilon
_{i}$. Moving to a subsequence if necessary, $d_{r}(x_{n_{k}},x_{n_{k+1}})<%
\epsilon _{n_{k}}$, for $n_{k}>N_{n_{k}}$. Thus, for $n_{k}>N_{n_{k}}$ we
have%
\begin{equation}
d_{r}(x_{n_{k}},\overline{x})\leq \sum_{i=k}^{\infty
}d_{r}(x_{n_{i}},x_{n_{i+1}})\leq \sum_{i=k}^{\infty }\epsilon _{i}
\end{equation}%
Hence, $\lim_{n_{k}\rightarrow \infty }d_{r}(x_{n_{k}},\overline{x})=0$ and,
since $\{x_{n_{k}}\}$ is a subsequence of a Cauchy sequence, $%
\lim_{n\rightarrow \infty }d_{r}(x_{n},\overline{x})=0$, as desired. The
result follows since $F$ is a contraction defined on the complete metric
space, $(X,d_{r})$.
\end{proof}

\bigskip

We next turn our attention to the case where $F$ is a set-valued pointwise
contraction. For $x$, $y\in X$ we will again let $d_{r}:X\times X\rightarrow 
%TCIMACRO{\U{211d} }%
%BeginExpansion
\mathbb{R}
%EndExpansion
_{+}$ be defined by

\begin{equation}
d_{r}(x,y):=\inf \{l(p):p\text{ is a rectifiable path from }x\text{ to }y\}
\end{equation}%
As discussed earlier, $d_{r}$ is a metric on $X$ and $(X,d_{r})$ is
complete, since $(X,d)$ is complete. Note that since $d_{r}$ is not
necessarily topologically equivalent to $d$, the space $(X,d_{r})$ might not
be compact. In what follows, we let $K(X):=\{Y\subseteq X:Y\neq \varnothing $%
, $Y$ compact$\}$ with respect to the metric $d_{r}$, and we again let $%
H_{r} $ denote the Hausdorff distance associated with $d_{r}$. We wish to
establish the following.

\begin{proposition}
Let $X$ be compact and rectifiably path-connected and let $F:X\rightarrow
K(X)$ be a set-valued pointwise contraction. Then \emph{Fix}$(F)\neq
\varnothing $.
\end{proposition}

\bigskip

The proof follows the steps in Ciesielski and Jasinski (2016b). Let $X$ be
as stated above. For any $Y\subset X$, $F:Y\rightarrow K(Y)$, and $x\in Y$,
define $\displaystyle s(F(x)):=\lim \sup_{y\rightarrow x}\frac{H(F(x),F(y))}{%
d(x,y)}$, if $x$ is a limit point of $Y$, and $s(F(x))=0$, otherwise.
Clearly, $F:X\rightarrow K(X)$ is a set-valued pointwise contraction if and
only if $s(F)<1$, for all $x\in X$. We have the following analog of Lemma 13.

\begin{lemma}
Let $Y$ be the range of a rectifiable path $p$ in $X$ and let $%
F:X\rightarrow K(X)$ be such that $s(F(x))<\beta \in \lbrack 0,1]$, for all $%
x\in Y$. Then $l(F\circ p)\leq \beta l(p)<\infty $, where 
\begin{equation*}
l(F\circ p):=\sup_{\pi \in \Pi }\dsum\limits_{i=1}^{n}H(F\circ
p(t_{i-1}),F\circ p(t_{i}))
\end{equation*}
\end{lemma}

\begin{proof}
As before, we let $l(p^{[a,b]})$ denote the restriction of a path $%
p:[0,1]\rightarrow X$ on $[a,b]$, where $0\leq a<b\leq 1$. We first show
that for any $\epsilon >0$ we have that $H(F\circ p(a),F\circ p(b))\leq
(\beta +\epsilon )l(p^{[a,b]})$.

Since $s(F(p(v)))\leq \beta <1$, for every $v\in \lbrack a,b]$, there exists 
$\delta _{v}>0$ and $B_{\delta _{v}}(v)$ such that for all $w\in B_{\delta
_{v}}(v)\cap \lbrack a,b]$ we have%
\begin{equation}
H(F\circ p(v),F\circ p(w))\leq (\beta +\epsilon )d(p(v),p(w))\ 
\end{equation}%
Since $[a,b]$ is compact, it has a minimal finite subcover $\mathcal{B}%
=\{B(v_{i}):v_{i}\in \lbrack a,b]\}$. Let $(v_{1},v_{3},...,v_{2n-1})$ list
the $v_{i}\in \lbrack a,b]$ in increasing order. Since $\mathcal{B}$ is
minimal, for all $1\leq i\leq n$ there exists $v_{2i}\in B(v_{2i-1})\cap
B(v_{2i+1})\cap (v_{2i-1},v_{2i+1})$, with $v_{0}=a\in B(v_{1})$ and $%
v_{2n}=b\in B(v_{2n-1})$. We thus obtain a partition $a=v_{0}\leq
v_{1}<v_{2}<...<v_{2n-1}\leq v_{2n}=b$ with $v_{2i},v_{2i+2}\in B(v_{2i+1})$%
, for all $0\leq i\leq n$. Inequality (24) then implies%
\begin{gather}
H(F\circ p(a),F\circ p(b))\leq \dsum\limits_{i=0}^{2n}H(F\circ
p(v_{i}),F\circ p(v_{i+1}))  \notag \\
\leq \dsum\limits_{i=0}^{2n}(\beta +\epsilon )d(p(v_{i}),p(v_{i+1}))\ \leq
(\beta +\epsilon )l(p^{[a,b]})
\end{gather}%
Next, let $\pi =\{t_{i}:t_{0}=0,...t_{i-1},t_{i},...t_{n}=1\}$ and $\epsilon
>0$ such that 
\begin{equation}
l(F\circ p)\leq \dsum\limits_{i=1}^{n}H(F\circ p(t_{i}),F\circ
p(t_{i+1}))+\epsilon
\end{equation}%
Inequality (25) implies 
\begin{eqnarray}
l(F\circ p) &\leq &\dsum\limits_{i=1}^{n}H(F\circ p(t_{i}),F\circ
p(t_{i+1}))+\epsilon  \notag \\
&\leq &\dsum\limits_{i=1}^{n}(\beta +\epsilon
)l(p^{[t_{i-1},t_{i}]})+\epsilon =(\beta +\epsilon )l(p)+\epsilon
\end{eqnarray}%
Since $\epsilon >0$ was arbitrary, the conclusion that $l(F\circ p)\leq
\beta l(p)$ follows.
\end{proof}

\bigskip 

We will also make use of the following.

\begin{lemma}
Let $X$ be compact and let $x,y\in X$. Suppose there is a rectifiable path $%
p:[0,1]\rightarrow X$. Then: (a) $p$ admits a reparametrization; i.e., a
path $\overline{p}:[0,1]\rightarrow X$ with the same length and range, such
that, f\textit{or all }$t\in \lbrack 0,1]$\textit{, }$l(\overline{p}%
^{[0,t]})=tl(\overline{p})$; (b) there is a geodesic; i.e., a rectifiable
path of minimum length, $L<\infty $, from $x$ to $y$.
\end{lemma}

\begin{proof}
\textit{(a)} Consider an arbitrary rectifiable path $p:[0,1]\rightarrow X$.
Rescale $t$ by letting $\displaystyle\overline{t}:=\frac{l(p^{[0,t]})}{l(p)}$
and define the mapping 
\begin{equation}
\overline{p}:=\left\{ \left( \frac{l(p^{[0,t]})}{l(p)},p(t)\right) :t\in
\lbrack 0,1]\right\}
\end{equation}%
We then have%
\begin{equation}
l(\overline{p}^{[0,\overline{t}]})=l(p^{[0,t]})=\overline{t}l(p)=\overline{t}%
l(\overline{p})
\end{equation}%
Thus, by rescaling the length of the segments of $p$, if needed, we can
always create a path, $\overline{p}$, with the same range and total length
such that $l(\overline{p}^{[0,t]})=tl(\overline{p})$, $t\in \lbrack 0,1]$.

\textit{(b)} Let $x,y\in X$ and let $L=\inf \{l(p):p$ is a rectifiable path
from $x$ to $y\}$. Choose a sequence of rectifiable paths $\{p_{n}\}$ from $%
x $ to $y$ such that $\lim_{n\rightarrow \infty }l(p_{n})=L$. By part 
\textit{(a)}, we can assume without loss of generality that, for all $n$, $%
l(p_{n})=tl(p_{n}^{[0,t]})$, for $t\in \lbrack 0,1]$. Since $X$ is compact,
applying Fact 1 to $\{p_{n}\}$, there is a subsequence $\{p_{n_{k}}\}$%
\textit{\ }converging uniformly to a rectifiable path $p:[0,1]\rightarrow X$
with $l(p)=L$.
\end{proof}

\bigskip

Next, we have the following.

\begin{definition}
A set-valued map $F:X\rightarrow K(X)$ is \textbf{shrinking }if for all $%
x,y\in X$ , $H(F(x),F(y))<d(x,y)$.
\end{definition}

The following Lemma establishes a connection between set-valued pointwise
contractions and set-valued maps satisfying the above shrinking property.

\begin{lemma}
Let $X$ be compact and let $F:(X,d)\rightarrow (K(X),H)$ be a set-valued
pointwise contraction. Then $F:(X,d_{r})\rightarrow (K(X),H_{r})$ is
shrinking.
\end{lemma}

\begin{proof}
Fix $x,y\in X$. We need to show that $H_{r}(F(x),F(y))<d_{r}(x,y)$. By Lemma
26(b), any two points in a compact metric space that can be joined by a
rectifiable curve, can be joined by a length minimizer. Thus, there is a
path $p:[0,1]\rightarrow X$ such that $p(0)=x$, $p(1)=y$, and $%
l(p)=d_{r}(x,y)=\inf \{l(p):p$ is a rectifiable path from $x$ to $y\}$.
Since $H_{r}(F(x),F(y))\leq l(F\circ p)$, it is sufficient to show that $%
l(F\circ p)<l(p)$.

Let $Y=p([0,1])\subseteq X$. Given $n>0$ and $\delta >0$, define the set 
\begin{equation}
M(n):=\{x\in Y:\frac{H(F(x),F(x^{\prime }))}{d(x,x^{\prime })}\leq \frac{n-1%
}{n}\text{, for all }x^{\prime }\in B_{1/n}(x)\subset Y\}
\end{equation}%
The set $M(n)$ is closed in $X$. Since $F$ is a pointwise set-valued
contraction, $Y=\cup _{n}M(n)$. Thus, $Y$ is an at most countable union of
closed sets. By the weak form of the Baire Category Theorem (see Royden and
Fitzpatrick, 2010), there exists $\overline{n}$ such that $int(M(\overline{n}%
))\neq \varnothing $ in $Y$. Since $p$ is continuous, there is an interval $%
[a,b]\subset p^{-1}\left( int(M(\overline{n}))\right) \subseteq \lbrack 0,1]$%
. For any $x\in Y=p([a,b])$ we have $\displaystyle s(F(x))\leq \frac{n-1}{n}$%
. Hence, by Lemma 16, $\displaystyle l(F\circ p^{[a,b]})\leq \frac{n-1}{n}%
l(p^{[a,b]})$. Moreover, $s(F(x))<1$, for all $x\in X$, so Lemma 16 also
implies that $l(F\circ p^{[v_{1},v_{2}]})\leq l(p^{[v_{1},v_{2}]})$, for any 
$0\leq v_{1}\leq v_{2}\leq 1$. Thus,%
\begin{eqnarray}
l(p) &=&l(p^{[0,a]})+l(p^{[a,b]})+l(p^{[b,1]})  \notag \\
&>&l(p^{[0,a]})+\frac{n-1}{n}l(p^{[a,b]})+l(p^{[b,1]})  \notag \\
&\geq &l(F\circ p^{[0,a]})+l(F\circ p^{[a,b]})+l(F\circ p^{[b,1]})  \notag \\
&=&l(F\circ p)
\end{eqnarray}%
We conclude that%
\begin{equation*}
H_{r}(F(x),F(y))\leq l(F\circ p)<l(p)=d_{r}(x,y)
\end{equation*}%
Therefore, $F$ is shrinking.
\end{proof}

\bigskip 

To complete the proof, we will make use of the following (Myers, 1945).\ 
\textit{Let }$X$\textit{\ be a compact metric space and for all }$n$\textit{%
, let }$p_{n}:[0,1]\rightarrow X$\textit{\ be a rectifiable path such that,
for all }$t\in \lbrack 0,1]$\textit{, }$l(p_{n}^{[0,t]})=tl(p_{n})$\textit{.
If }$L:=\lim \inf_{n\rightarrow \infty }l(p_{n})<\infty $\textit{, then
there is a subsequence }$\{p_{n_{k}}\}$\textit{\ converging uniformly to a
rectifiable path }$p:[0,1]\rightarrow X$\textit{\ with }$l(p)\leq L$\textit{.%
}

\bigskip 

\begin{proof}
Returning to the proof of the Proposition, let 
\begin{equation}
L=\inf \{H_{r}(x,F(x)):x\in X\}
\end{equation}%
It remains to show that there exists an $x^{\ast }\in X$ such that $%
H_{r}(x^{\ast },F(x^{\ast }))=L=0$; i.e., $x^{\ast }\in F(x^{\ast })$. By
Lemma 19, $H_{r}(F(x),F(y))<d_{r}(x,y)$, for all $x,y\in X$. Fix a sequence $%
\{x_{n}\}$ such that, for all $n$, $x_{n+1}\in F(x_{n})$ is chosen such that 
$d_{r}(x_{n},x_{n+1})\leq H_{r}(x_{n},F(x_{n}))$. Then, 
\begin{equation*}
\lim \inf_{n\rightarrow \infty }d_{r}(x_{n},x_{n+1})\leq \lim
\inf_{n\rightarrow \infty }H_{r}(x_{n},F(x_{n}))=L
\end{equation*}%
By rectifiability and Lemma 17(b), for every $n$, there exists a path $%
p_{n}:[0,1]\rightarrow X$, $x_{n}\longmapsto x_{n+1}\in F(x_{n})$, of
(minimizing) length $d_{r}(x_{n},x_{n+1})=\inf \{l(p_{n}):p_{n}$ is a
rectifiable path from $x_{n}$ to $x_{n+1}\}$. Thus, 
\begin{equation*}
\lim \inf_{n\rightarrow \infty }l(p_{n})=\lim \inf_{n\rightarrow \infty
}d_{r}(x_{n},x_{n+1})\leq L
\end{equation*}%
Using Lemma 17(a) to reparametrize, if necessary, we have that $%
l(p_{n})=tl(p_{n}^{[0,t]})$, for all $n$ and $t\in \lbrack 0,1]$. Myers
(1945) implies that $\{p_{n}\}$ admits a subsequence $\{p_{n_{k}}\}$ that
converges uniformly to a rectifiable path $p:[0,1]\rightarrow X$ with $%
l(p)\leq L$. Let $x^{\ast }=\lim_{k\rightarrow \infty
}p_{n_{k}}(0)=\lim_{k\rightarrow \infty }x_{n_{k}}$. Then $p$ is a path from 
$p(0)=x^{\ast }$ to $p(1)=\lim_{k\rightarrow \infty
}p_{n_{k}}(1)=\lim_{k\rightarrow \infty }x_{n_{k+1}}\in \lim_{k\rightarrow
\infty }F(x_{n_{k}})=F(x^{\ast })$, where the last implication follows from
the fact that, since it is compact-valued and Lipschitz, $F$ has a closed
graph. Thus, $H_{r}(x^{\ast },F(x^{\ast }))\leq l(p)\leq L$. Finally, note
that we must have $L=0$, since otherwise $H_{r}(F(x^{\ast }),F(F(x^{\ast
})))<H_{r}(x^{\ast },F(x^{\ast }))$, contradicting that $L=\inf
\{H_{r}(x,F(x)):x\in X\}$. Hence, $H_{r}(x^{\ast },F(x^{\ast }))=0$, and $%
x^{\ast }\in F(x^{\ast })$ as desired. We conclude that \emph{Fix}$(F)\neq
\varnothing $.
\end{proof}

\section{An Application: Pure-Strategy Nash Equilibria in Metric Spaces}

A game in strategic form, $\Gamma $, is defined by $\Gamma =\left\langle
I,X_{i},u^{i}\right\rangle $, where $I=\{1,...,n\}$ is the set of players, $%
X_{i}$ is player $i$'s (pure) strategy set, and $u_{i}:X_{1}\times ...\times
X_{n}\rightarrow 
%TCIMACRO{\U{211d} }%
%BeginExpansion
\mathbb{R}
%EndExpansion
$ is the payoff of player $i$, $i=1,...,n$. Each player's payoff depends on
his strategy, as well as on that of the other players. We assume that the
players' strategy sets are metric spaces, $(X_{i},d_{i})$, and use $\mathbf{x%
}=(x_{1},...,x_{n}):=(x_{i},\mathbf{x}_{-i})$ to denote the strategy
profile, where $x_{i}\in X_{i}$, and $\mathbf{x}_{-i}\in X_{1}\times
...\times X_{i-1}\times X_{i+1}\times ...\times X_{n}$, for all $i$. Thus,
the payoff of each player $i$ can be written as $(x_{i},\mathbf{x}%
_{-i})\mapsto u_{i}(x_{i},\mathbf{x}_{-i})$. Let $X=X_{1}\times ...\times
X_{n}$, and let $d$ be the product metric associated with the metrics $d_{i}$%
, $i=1,...,n$, and similarly for $d_{r}$. For all $i$, the \textbf{best
response of player }$i$, $BR_{i}:X\rightarrow X_{i}$, is defined by the
set-valued map

\begin{equation}
BR_{i}(\mathbf{x}):=\left\{ x_{i}\in X_{i}:u_{i}(x_{i},\mathbf{x}_{-i})\geq
u_{i}(x_{i}^{\prime },\mathbf{x}_{-i})\text{, for all }x_{i}^{\prime }\in
X_{i}\right\} 
\end{equation}%
The \textbf{best response for the game }$\Gamma $, $BR:X\rightarrow X$, is
defined by the set-valued map $BR(\mathbf{x}):=BR_{1}(\mathbf{x})\times
...\times BR_{n}(\mathbf{x})$. We have the following. A \textbf{(pure
strategy) Nash equilibrium }for $\Gamma $ is a strategy profile $\mathbf{x}%
^{\ast }\in X$ such that $\mathbf{x}^{\ast }\in BR(\mathbf{x}^{\ast })$.

The finite Cartesian product of compact sets is compact, and similarly for
rectifiably path-connected sets. The following follows directly from the
existence results established in the previous sections.

\begin{corollary}
Suppose that either one of the following conditions hold: (a) $(X,d)$ is
compact and $r$-chainable, and $BR(\mathbf{x})$ is a local set-valued
contraction on $X$; (b) $(X,d)$ is complete and rectifiably path connected,
and $BR(\mathbf{x})$ is a uniform pointwise set-valued contraction on $X$ ;
(c) $(X,d_{r})$ is compact and rectifiably path connected, and $BR(\mathbf{x}%
)$ is a pointwise set-valued contraction on $X$. Then \emph{Fix}$(BR(\mathbf{%
x}))\neq \varnothing $.
\end{corollary}

\bigskip

\noindent \textbf{Declarations}

No funds, grants, or other support was received. The authors have no
conflicts of interest to declare that are relevant to the content of this
article.

\end{document}